\documentclass{article}
\usepackage{latexsym,amssymb,amsmath}
\usepackage{amsfonts}

\numberwithin{equation}{section}
\newtheorem{theorem}{Theorem}[section]

\newtheorem{proposition}[theorem]{Proposition}

\newtheorem{definition}[theorem]{Definition}
\newtheorem{remark}[theorem]{Remark}
\newenvironment{proof}[1][Proof]{\textbf{#1.} }{\ \mathbb Rule{0.5em}{0.5em}}

\newcommand{\po}{{\mathbb P}}
\newcommand{\Ff}{{\mathbb F}}
\newcommand{\om}{{\omega}}

\begin{document}

\title{A note on a result of Liptser-Shiryaev}
\author{Benedetta Ferrario\\
 Dipartimento di Matematica - Universit\`a di Pavia}
\date{\today}
\maketitle

\begin{abstract}
Given two stochastic equations with different drift terms, 
under very weak assumptions 
Liptser and Shiryaev provide the equivalence of the laws of the
solutions to these equations   by means of
Girsanov transform. Their assumptions involve both the drift
terms. We are interested in the same result but with the main 
assumption involving only the difference of the drift terms. Applications of
our result will be presented in the finite  as well
as in the infinite dimensional setting.
\end{abstract}
\maketitle

\section{Introduction}\label{intro}
Let us consider the It\^o equation
\begin{equation}\label{nlin}
 dX(t)=b(t,X(t))\  dt+\sigma(t,X(t))\ dW(t), 
 \quad  X(0)=x.
\end{equation}
If we know that there exists a solution, we ask about uniqueness and
characterization of its law. 
We can look at equation \eqref{nlin} as a modification of equation
\begin{equation}\label{lin}
dZ(t)=a(t,Z(t))\ dt+\sigma(t,Z(t))\ dW(t), \quad  Z(0)=x.
\end{equation}
by a change of the drift term.
Equation \eqref{lin} is a ''good'' reference equation, for
which existence and uniqueness hold true.
Since  these two equations differ only in the drift terms,
a  classical tool to study equation
\eqref{nlin} is the Girsanov transform. 

In \cite{LS}, Chapter 7 is
devoted to this problem, where Liptser and Shiryaev 
investigate the relation between the laws of processes solving
equations \eqref{nlin} and \eqref{lin}. 
In this paper,
we address the same problem.

As far as the results are concerned, first in dimension one
 we prove  results similar to 
\cite{LS} but with different 
assumptions; in fact, our hypotheses involve the difference $b-a$
whereas in \cite{LS} involve separately $a$ and $b$.
Then, we consider the case of dimension bigger than one.
Our analysis includes the 
uniqueness problem, not tackled in \cite{LS}.
Moreover we extend these results 
to the infinite dimensional setting, whereas \cite{LS} 
deals only with the finite dimensional case.
Here, when 
we say finite dimensional we mean that the state space is finite 
dimensional, i.e. the
unknown $X$
 is a vector process with a finite number ($d<\infty$) of components; 
this models stochastic differential equations
on the state space $\mathbb R^d$. 
However, the infinite dimensional setting is related to abstract
models of  
stochastic partial differential equations
 (see, e.g., the book by Da Prato and Zabczyk \cite{dpz}).
Actually, the infinite dimensional setting is one of the main motivations of
 our study, as it will be explained in Section \ref{sAPP}.

As far as the techniques are concerned, in some parts 
our proofs are shorter than in \cite{LS}, in the sense that even with
the same assumption of \cite{LS} we get the results of \cite{LS}  with
shorter proofs.

Now, we explain how the paper is organized.
We start our exposition with the  one dimensional setting. Extension
to dimension bigger than one is in Section \ref{bigd}.
After the basic results presented in Sections \ref{sPRE} and \ref{easy}, 
 we shall analyze uniqueness in law  in Section \ref{coro},
 the absolute continuity  in Section \ref{mains}  
 and the equivalence of the laws in Section \ref{s-equi}. 
In Section \ref{conclu} our results will be compared with those in \cite{LS}.
In the  final
section the novelty 
of our results will be discussed, also in the infinite dimensional setting.

\section{Preliminaries}\label{sPRE}
We set our problem as in the book of Liptser and Shiryaev \cite{LS}, 
that is in a setting more general than 
\eqref{nlin}-\eqref{lin}.

Let $(\Omega, \Ff,\po)$ be a probability space and  $\{\Ff_t\}_{t\ge
  0}$ a filtration. We will always assume that 
the probability space is complete and the filtration is right
continuous. 
We denote by $\mathbb E$ the expectation with respect to the measure
$\po$, and by
 $\Ff_T(X)$ the $\sigma$-algebra generated by $\{X(u), 0\le u \le
T\}$.

When dealing with a Polish space,
i.e. a complete separable metric space, 
the $\sigma$-algebra associated is the Borel $\sigma$-algebra.
In particular,  for $0<t\le T$ let 
$\mathcal B_t$ be the $\sigma$-algebra of Borelian subsets of
$C([0,t];\mathbb R)$.
We say that a measurable functional 
$\phi:[0,T]\times C([0,T];\mathbb R)\to \mathbb R$ 
is  non anticipative
 if, for each $t \in [0,T]$, $\phi(t,\cdot)$  is $\mathcal B_{t}$-measurable.

The two equations to deal with are
\begin{equation}\label{equX}
 dX(t)=  b(t,X)\ dt+ \sigma(t,X)\ dW(t); \qquad X(0)=x 
\end{equation}
\begin{equation}\label{equZ}
 dZ(t)=  a(t,Z)\ dt+ \sigma(t,Z)\ dW(t); \qquad Z(0)=x 
\end{equation}

Here, $a$, $b$ and $\sigma$ are  non 
anticipative measurable functionals.
$W$ is a  Wiener process with respect to the
stochastic basis $(\Omega, \Ff,\{\Ff_t\}, \po)$.

We need to recall what is a weak or strong solution.
We consider  processes $X$ with a.e. path in $C([0,T];\mathbb R)$, which are
adapted to the 
filtration $\{\Ff_t\}_{t\ge 0}$ and solve  equation \eqref{equX} a.s.:
\begin{equation}
 X(t)=x+ \int_0^t  b(s,X)\ ds+\int_0^t\sigma(s,X)\ dW(s)  \qquad \po-a.s.
\end{equation}
for every $t \in [0,T]$.
It is necessary that 
$$
 \po\{\textstyle\int_0^T|b(s,X)|ds<\infty\}
 =\po\{\int_0^T \sigma(s,X)^2 ds<\infty\}
 =1. 
$$
For simplicity, we fix the initial data $x \in \mathbb R$; however,
our results can be extended to cover the case of random initial data.
\begin{definition}[weak solution]
We say that there exists  a weak solution to equation \eqref{equX} 
if there exist a stochastic basis 
$(\Omega, \Ff,\{\Ff_t\}_{t\ge 0},\po)$,
an  $\{\Ff_t\}$-Wiener process $W$
and 
 an $\{\Ff_t\}$-adapted process $X$ defined in it   such that
$X$ solves  equation \eqref{equX} $\po$-a.s.\\
We denote this solution by the triplet 
$\big(X,\ (\Omega, \Ff,\{\Ff_t\},\po), \ W \big)$.
\end{definition}
On the other hand, if $X$ 
solves \eqref{equX} on a (a priori) given stochastic basis 
$(\Omega, \Ff,\{\Ff_t\}_{t\ge 0},\po)$
with a given Wiener process $W$,
we have a strong solution. Therefore the Wiener 
process and the filtration are not part of the
solution itself but are assigned.

\begin{definition}[strong solution] 
We say that there exists  a strong solution 
to equation \eqref{equX} 
if, given any stochastic basis 
$(\Omega, \Ff,\{\Ff_t\}_{t\ge 0},\po)$ and
$\{\Ff_t\}$-Wiener process $W$, 
there exists an $\{\Ff_t\}$-adapted process $X$   such that
$X$ solves  equation \eqref{equX} $\po$-a.s. 
\end{definition}

Moreover, we have two kinds of uniqueness. 

\begin{definition}[uniqueness in law] 
We say that uniqueness in law holds for equation \eqref{equX} if 
any two processes solving equation \eqref{equX}
with the same initial data have the same law.
\end{definition}

\begin{definition}[pathwise uniqueness] $\ $
We say that pathwise  uniqueness  holds for equation \eqref{equX} 
if  given two processes $X$ and $X^\prime$ solving equation \eqref{equX}
with the same initial data and defined with respect to the same
stochastic basis $(\Omega,\Ff,\{\Ff_t\}_{t\ge 0},\po)$ and 
Wiener process, 
we have 
$\po\{X(t)= X^\prime(t) \text{ for all } t\}=1$.
\end{definition}

In the following we shall assume that equation \eqref{equZ} has a
unique strong solution; uniqueness has to be understood as pathwise
uniqueness.
 But, 
a result of Cherny (see \cite{che}) says that
 uniqueness in law,  together with the strong existence, guarantees
the pathwise uniqueness.
Hence, we could simply assume existence of a strong solution and
uniqueness in law.
\\
On the other hand, from now on saying uniqueness of a  weak
solution we will mean 
uniqueness in law, unless otherwise specified.

Therefore,  the coefficients 
$a$ and $\sigma$ are required to 
satisfy the usual growth and Lipschitz conditions (see, e.g.,
\cite{LS}), that is
\[
{\mathbf{ [A1]}}
\left[
\begin{array}{l}
\; \exists \text{ constants } L_1,L_2 \text{ and a function }
K \text{ non decreasing and right continuous,}\\
\text{ with }0 \le K(s)\le 1, 
\text{ such that   }
\\[2mm]
 a(t,Y)^2+\sigma(t,Y)^2 
\le L_1 \int_0^t [1+Y(s)^2] dK(s)+L_2[1+ Y(t)^2 ]
  \\\hspace*{66mm} \forall t \in [0,T], Y \in C([0,T];\mathbb R)
\\\text{ and } \\
 |a(t,Y_1)-a(t,Y_2)|^2+|\sigma(t,Y_1)-\sigma(t,Y_2)|^2
\\
\hspace*{30mm}
\le L_1\int_0^t |Y_1(s)-Y_2(s)|^2dK(s)+L_2 |Y_1(t)-Y_2(t)|^2 
 \\\hspace*{66mm} \forall t \in [0,T], Y_1,Y_2 \in C([0,T];\mathbb R)
\end{array}
\right.
\]
Moreover,
the coefficients $a,b$ and $\sigma$ are such that
\[
{\mathbf{ [A2]}}
\left[
\begin{array}{l}
\; \exists \text{ a measurable functional } \gamma \text{ which is 
non anticipative  finite and  such that }\\ 
\quad\sigma(s,Y)\gamma(s,Y)=b(s,Y)-a(s,Y)\qquad
\forall s\in [0,T], Y \in C([0,T];\mathbb R).
\end{array}
\right.
\]

Few technical details:
from now on, we consider only finite time intervals $[0,T]$. 
Then, the law of a process 
solving equation \eqref{equX} or \eqref{equZ} is a probability measure on
$\mathcal B_T$.
Moreover,
if there is uniqueness in law for an equation with drift term $a$ we
denote by $\mu^a$ this unique law (unless otherwise stated).
If a measure $\nu_1$ is absolutely continuous with respect to a measure $\nu_2$ we write
$\nu_1\prec\nu_2$; if they are equivalent, i.e. $\nu_1\prec\nu_2$ and $\nu_2\prec\nu_1$, we 
write $\nu_1\sim\nu_2$.

\section{An easy case}\label{easy}
In this section, we prove a result of equivalence of laws for
equations
\eqref{equZ} and \eqref{equX}  but in the particular case 
of $b=a+g$ with a strong assumption on $\sigma$ and $g$. The proof is
based on classical tools of Girsanov transform and Novikov condition.

Instead of equation \eqref{equX}, let us  consider
\begin{equation}\label{equY}
 dY(t)=  a(t,Y)\ dt+ g(t,Y)\ dt\ +\sigma(t,Y)\ dW(t), \qquad Y(0)=x,
\end{equation}
where $g$ is a non anticipative measurable functional. Moreover, we assume that
there exists a 
finite and    non anticipative measurable functional
$\alpha$  such that
\begin{equation}\label{sistema}
 \sigma(s,Y)\alpha(s,Y)=g(s,Y)
\end{equation}
for each $s\in [0,T]$ and $Y \in C([0,T];\mathbb R)$. 

\begin{remark}\label{oss-su-s+}
Relationship \eqref{sistema}
is a compatibility condition; 
it means that when $\sigma$ vanishes, also $g$ must
vanish. 
In this case, $\alpha$ may be chosen arbitrarily in order to
satisfy \eqref{sistema}.
But, as we shall see, 
 in the Girsanov transform $\alpha$ takes into account the change of
drift between equations  \eqref{equY} and \eqref{equZ}. Therefore we
are interested only in the solution $\alpha$ of \eqref{sistema} which
vanishes when $g=0$, i.e. when the two drift terms are the same. 
Hence, from now on we consider
\begin{equation}\label{laAlfa}
\alpha(s,Y)=\sigma^+(s,Y)g(s,Y)
\end{equation}
where
$$
\sigma^+(s,Y)=\begin{cases}\dfrac 1{\sigma(s,Y)},& \sigma(s,Y)\neq 0\\  
                           0,& \sigma(s,Y)=0  \end{cases}
$$
\end{remark}

We have the following result.
\begin{theorem}\label{semplice}
Assume 
 there exists a unique weak solution 
$\big(Z,\ (\Omega, \Ff,\{\Ff_t\},\po), \ W \big)$
to equation \eqref{equZ}.
If 
\begin{equation}\label{unifo}
 \sup_{X\in C([0,T];\mathbb R)}\int_0^T \alpha(s,X)^2ds=c<\infty,
\end{equation}
then  equation \eqref{equY} has a weak solution, which is unique in
law.
Moreover, the law of the process $Z$ is equivalent to the law of the
process solving \eqref{equY}, that is $\mu^a \sim \mu^{a+g}$.
In particular
\begin{equation}\label{RNperY}
\frac{d\mu^{a+g}}{d\mu^{a}\;\;}(Z)=
\mathbb E\Big[e^{\textstyle\int_0^T \alpha(s,Z)dW(s)
-\frac 12 \textstyle\int_0^T \alpha(s,Z)^2ds}\big|\Ff_T(Z)\Big]
\end{equation}
$\po$-a.s.
\end{theorem}
\proof
Because of \eqref{unifo} we have that 
$$
  \mathbb E\Big[e^{\textstyle\frac 12 \int_0^T \alpha(s,Z)^2 ds}\Big]
   \le e^{\textstyle\frac c2}<\infty.
$$
This is Novikov condition, which allows to apply Girsanov
transform. More precisely (see \cite{KS}), Novikov condition makes sure that
the process $\delta$ defined by 
$$
 \delta_t=e^{\textstyle\int_0^t \alpha(s,Z)dW(s)  
-\frac 12 \textstyle\int_0^t \alpha(s,Z)^2 ds}, \; 0\le t\le T,
$$
is a martingale. To highlight the dependence on $Z$ and $W$ we will
often write $\delta_T$ as $\delta_T(Z,W)$.
We define a new probability measure on 
$(\Omega, \mathcal F_T)$ by  $d\po^*= \delta_T(Z,W) d\po$. 
Then Girsanov theorem (see \cite{gir}) tells us that
$$
W^{*}(t)=W(t)-\int_0^t \alpha(s,Z) \ ds \ , \qquad t \in [0,T],
$$
is a  Wiener process with respect to
$(\Omega, \Ff,\{\Ff_t\},\po^{*})$; 
substituting into equation \eqref{equZ} we get
$$
 Z(t)=x +\int_0^t a(s,Z)\ ds + \int_0^t g(s,Z)\  ds
  +\int_0^t\sigma(s,Z)dW^{*}(s).
$$
This means that 
$\big(Z,\ (\Omega, \Ff,\{\Ff_t\},\po^*), \ W^* \big)$
is a weak solution of equation \eqref{equY}.

For any Borelian subset $\Lambda$ of $C([0,T];\mathbb R)$, 
set $\mathcal L_Y(\Lambda)=\po^*\{Z \in \Lambda\} $ and 
$\mu^{a}(\Lambda)=\po\{Z \in \Lambda\}$. Then
$\mathcal L_Y \prec \mu^{a}$, since
$\po^*\prec \po$ by construction. Moreover, 
consider the  random variable $ \mathbb E[\delta_T(Z,W)|\Ff_T(Z)]$; it
is $\Ff_T(Z)$-measurable and therefore there exists a  
$\mathcal B_T$-measurable non negative function 
$D: C([0,T];\mathbb R)\to \mathbb R$ such that
$D(Z(\omega))=\mathbb E[\delta_T(Z,W)|\Ff_T(Z)](\omega)$ for $\po$-a.e. $\omega$.
Now, we have
\begin{multline}
 \mathcal L_Y(\Lambda)=\po^*\{Z \in \Lambda\}
 =\int\limits_{\{Z \in \Lambda\}}\delta_T(Z,W)\ d\po
 = \int\limits_{\{Z \in \Lambda\}}\mathbb E\big[\delta_T(Z,W)|\Ff_T(Z)\big]\ d\po
\\
 =\int\limits_{\{Z \in \Lambda\}}D(Z) \ d\po
 =\int_\Lambda D(z)\ d\mu^a(z)
\end{multline}
Hence 
$$
 \frac{d\mathcal L_Y}{d\mu^a}(Z)=D(Z) \quad \text{ for } Z \in
 C([0,T];\mathbb R).
$$
This proves \eqref{RNperY},
as soon as we have uniqueness in law for equation \eqref{equY}.

Viceversa, any weak solution 
$\big(Y,(\tilde \Omega, \tilde\Ff,\{\tilde \Ff_t\},\tilde\po),\tilde W\big)$
of equation \eqref{equY} 
gives rise to a weak solution 
$\big(Y,(\tilde \Omega, \tilde\Ff,\{\tilde \Ff_t\},\tilde\po^*),\tilde W^*\big)$
of equation \eqref{equZ}, with a similar
expression of the Radon-Nikodym derivative  (only a change of sign
appears). Indeed, thanks to \eqref{unifo}
\begin{equation}\label{p-den}
\Hat\delta_t(Y,\tilde W)=e^{\textstyle-\int_0^t \alpha(s,Y)d\tilde W(s)
-\frac 12 \textstyle\int_0^t \alpha(s,Y)^2 ds}
\end{equation}
is a martingale; define $d\tilde\po^*=\Hat\delta_T(Y,\tilde W) d\tilde \po$ and 
$ \tilde W^{*}(t)=\tilde W(t)+\int_0^t \alpha(s,Y) \ ds$.
Then, $\tilde W^{*}$ is a  Wiener process with respect to $\tilde \po^{*}$ and
\begin{equation}\label{altraRN}
 \mu^a(\Lambda)=\tilde\po^*\{Y\in \Lambda\}
 = \int\limits_{\{Y\in \Lambda\}} \Hat\delta_T(Y,\tilde W) \ d\tilde \po.
\end{equation}
Now, suppose there exist two different weak solutions of equation \eqref{equY}:
$$
 \big(Y_i,
 (\tilde \Omega_i, \tilde\Ff_i,\{\tilde \Ff_{{i}_t}\},\tilde\po_i),
  \tilde W_i\big)\qquad
 i=1,2.
$$
We have that $d\tilde\po_i^*=\Hat\delta_T(Y_i,\tilde W_i) d\tilde \po_i$;
moreover, 
\[\begin{split}
 \Hat\delta_t(Y_i,\tilde W_i)&=e^{\textstyle-\int_0^t \alpha(s,Y_i)d\tilde W_i(s)
-\frac 12 \textstyle\int_0^t \alpha(s,Y_i)^2 ds}
\\&=
e^{\textstyle-\int_0^t \alpha(s,Y_i)d\tilde W_i^*(s)
+\frac 12 \textstyle\int_0^t \alpha(s,Y_i)^2 ds}=:
\frac 1{\underline \delta_t(Y_i,\tilde W^*_i)}.
\end{split}\]
Again \eqref{unifo} provides that 
$\underline \delta_T(Y_i,\tilde W^*_i)$ is well defined.
Then, $d\tilde\po_i=\underline\delta_T(Y_i,\tilde W^*_i) d\tilde
\po^*_i$. Now, 
uniqueness in law for the solution of equation \eqref{equZ} 
means that the joint distribution of $(Y_1,W^*_1)$ is the same as of 
$(Y_2,W^*_2)$ (see \cite{che}  Th. 3.1).
Then, we get
\[\begin{split}
\tilde P_1(Y_1\in \Lambda)
  &=\int_{\tilde \Omega_1} \underline \delta_T(Y_1,\tilde W^*_1) 
  \mathbb I_{\{Y_1  \in \Lambda\}}d\tilde \po^*_1 
   \\&=
    \int_{\tilde \Omega_2} \underline \delta_T(Y_2,\tilde W^*_2) 
  \mathbb I_{\{Y_2  \in \Lambda\}}d\tilde \po^*_2 
  = 
  \tilde P_2(Y_2\in \Lambda)
\end{split}\]
for any Borelian subset $\Lambda$ of $ C([0,T];\mathbb R)$; 
here ${\mathbb I}_\cdot$ is the indicator  function. 
Thus, we have uniqueness in law for equation \eqref{equY}.
\hfill $\Box$

\begin{remark}\label{ledensi}
i) The expression  \eqref{RNperY} can be written as
$$
\frac{d\mu^{a+g}}{d\mu^{a}\;\;}(Z)=
\mathbb E\Big[e^{\textstyle\int_0^T \alpha(s,Z)dW(s)}\big|\Ff_T(Z)\Big]
e^{-\frac 12 \textstyle\int_0^T \alpha(s,Z)^2ds}.
$$
The same holds for other similar expressions of  Radon-Nikodym
derivatives appearing later on.
\\ii) Consider the assumptions of Theorem \ref{semplice}. Then, 
given a weak solution
$\big(Y,(\tilde \Omega, \tilde\Ff,\{\tilde \Ff_t\},\tilde\po),\tilde W\big)$  
of equation \eqref{equY}, from \eqref{altraRN} in the previous proof  we have
\begin{equation}\label{a:su:a+g}
 \frac{d\mu^{a}\;\;}{d\mu^{a+g}}(Y)=
 \tilde{ \mathbb E}\Big[e^{\textstyle-\int_0^T \alpha(s,Y)d\tilde W(s)
 -\frac 12 \textstyle\int_0^T \alpha(s,Y)^2 ds}\big|\Ff_T(Y)\Big]
\end{equation}
$\tilde \po$-a.s.
\end{remark}

\section{Uniqueness in law}\label{coro}
According to Remark \ref{oss-su-s+},
if {\bf [A2]} holds true we set
$$
 \gamma(s,Y)=\sigma^+(s,Y)[b(s,Y)-a(s,Y)].
$$

We have the following
\begin{proposition}\label{unici}
Assume {\bf [A1]} and {\bf [A2]}.\\
If there exist two weak solutions
$\big(X, (\Omega, \Ff,\{\Ff_t\},\po), \ W\big)$ and 
$\big(X^\prime, (\Omega^\prime, \Ff^\prime,\{\Ff^\prime_t\},\po^\prime), 
\ W^\prime\big)$ to equation \eqref{equX}, 
 with the same initial data $x$, 
such that
\begin{equation}\label{unic}
 \po\{\textstyle \int_0^T \gamma(s,X)^2 ds<\infty
   \}=
 \po^\prime\{\textstyle \int_0^T \gamma(s,X^\prime)^2 ds<\infty
 \}=1,
\end{equation}
 then the laws of $X$ and $X^\prime$ are the same. 
\end{proposition}
\proof 
Consider the first solution $\big(X, (\Omega, \Ff,\{\Ff_t\},\po), \ W\big)$.
According to {\bf [A1]} there exists a 
solution $Z$ of equation \eqref{equZ} with respect 
to the stochastic basis $(\Omega, \Ff,\{\Ff_t\},\po)$
and the Wiener process $W$.
For any integer $n\ge 1$, define the truncation function 
$$
 \chi^n_t(Z)= 
\begin{cases} 
1 &\text{ if }\int_0^t  \gamma(s,Z)^2ds < n,
\\
0 &\text{ otherwise. } 
\end{cases}
$$
We have that   
$$
 \sup_{Z \in C([0,T];\mathbb R)}
  \int_0^T \chi^n_s(Z) \gamma(s,Z)^2 ds\le n.
$$
We use Theorem \ref{semplice} with 
$g(s,Y)=\chi^n_s(Y)[b(s, Y)-a(s,Y)]$ in order to 
get that
the new  equation
\begin{multline}\label{eqN}
 Y(t)=x+ \int_0^t a(s,Y)\ ds+\int_0^t \chi^n_s(Y)[b(s,Y)-a(s,Y)]\  ds
\\ +\int_0^t\sigma(s,Y)dW(s)
\end{multline}
has a unique weak solution. For short, 
we denote its law by $\mu^{b,n}$ and we have
$\mu^{b,n}\prec \mu^a$, with
$\mu^{b,n}(\Lambda)=\po^{*n}\{Z\in \Lambda\}$,
$d\po^{*n}=\rho^n_T(Z,W) d\po$, 
and the martingale $\rho^n=\rho^n(Z,W)$ defined by
$\rho^n_t=e^{\textstyle\int_0^t \chi^n_s(Z) \gamma(s,Z)dW(s) 
-\frac 12 \textstyle\int_0^t \chi^n_s(Z) \gamma(s,Z)^2 ds}$. In particular,
\[
\mathbb E[\rho_t^n(Z,W)]=\mathbb E[\rho_0^n(Z,W)]=1 \qquad\text{ for all }t,
\]
and
\begin{equation}\label{den-n}
\frac{d\mu^{b,n}}{d\mu^{a}}(Z)=
\mathbb E\big[e^{\textstyle\int_0^T \chi_s^n(Z)\gamma(s,Z)dW(s)
-\frac 12 \textstyle\int_0^T\chi_s^n(Z)
\gamma(s,Z)^2 ds}|\Ff_T(Z)\big]
\end{equation}
$\po$-a.s..\\
This holds for any $n$ integer. 
Therefore we have uniquely defined the sequence $\{\mu^{b,n}\}_{n=1}^\infty$.

On the other hand, 
we can define
a process solving equation \eqref{eqN} {\it with the Wiener process} $W$.
Let us
define the sequence of stopping times (depending on the process $X$)
$$
 \tau^n= \inf\{t\in [0,T]: \chi^n_t(X)=0 \} \wedge T 
$$
considering the infimum to be $+\infty$ when the set is empty.\\
Given any $n$,  
$(\Omega, \Ff,\{\Ff_t\},\po)$ and  $ W$,
if $\big(X, (\Omega, \Ff,\{\Ff_t\},\po), \ W\big)$  
is a  weak solution to equation \eqref{equX} then
equation
\begin{multline}\label{Neq}
X^n(t)=X(t\wedge\tau^n)\\+\int_0^t[1-\chi^n_s(X)]a(s,X^n)\ ds
+\int_0^t[1-\chi^n_s(X)]\sigma(s,X^n)\ dW(s)
\end{multline}
has a unique strong solution $X^n$, thanks to
assumption {\bf [A1]}. 
Moreover, by It\^o calculus we get that this process $X^n$ solves  \eqref{eqN}.
Hence, $\mu^{b,n}$ coincides with the law of $X^n$.
From \eqref{Neq} we have ($\po$-a.s)
\[
X^n(t)=\begin{cases}X(t) \quad \text{ on } \{\tau^n\ge t\}\\
          X(\tau^n)+\int_{\tau^n}^t a(s,X^n)ds+\int_{\tau^n}^t \sigma(s,X^n)dW_s 
              \quad \text{ on }   \{\tau^n< t\}\end{cases}
\]
In particular, $ X=X^n$ on the set $\{\tau^n=T\}\supseteq \{\chi_T^n(X)=1\}$.
According to  \eqref{unic} we have 
$\displaystyle\lim_{n \to \infty}\po\{ \chi_T^n(X)=0\}=0$. 
Hence
\[
 \po\{ \| X-X^n\|_{C([0,T];\mathbb R)}>0\} \le \po\{ \chi_T^n(X)=0\} \to 0  \qquad
 \text{ as } n \to \infty. 
\]
Therefore, ${\mu^{b,n}}$ 
converges to the law of $X$ in the metric of total variation.

If we start from another solution 
$\big(X^\prime, (\Omega^\prime, \Ff^\prime,
\{\Ff^\prime_t\},\po^\prime), \ W^\prime\big)$ fulfilling \eqref{unic},
we would consider
the solution 
$\big(Z^\prime, (\Omega^\prime, \Ff^\prime,
\{\Ff^\prime_t\},\po^\prime), \ W^\prime\big)$
to equation \eqref{equZ}, giving  the same $\mu^{b,n}$. Indeed, there
is uniqueness in law for both equations \eqref{equZ} and \eqref{equY}.
Hence, $\mu^{b,n}$ converges to the law of $X^\prime$.
Since the limit  of $\mu^{b,n}$ is unique, we conclude that the laws of $X$
and $X^\prime$ are the same.
\hfill $\Box$

\begin{remark}\label{sol:deb}
Given a weak solution 
$\big(X, (\Omega, \Ff,\{\Ff_t\},\po), \ W\big)$
of equation \eqref{equX}, it is easier 
to construct a  solution 
$\big(X^n,\ (\Omega, \Ff,\{\Ff_t\},\po), \ W \big)$ 
of equation \eqref{eqN}
when the noise is independent of the unknown,
 i.e. $\sigma(t,X)=\sigma(t)$. 
Indeed, 
we look for a process solving
\[
 X^n(t)=x+\int_0^t a(s,X^n)\ ds+\int_0^t \chi^n_s(X^n)[b(s,X^n)-a(s,X^n)]\  ds
 +\int_0^t\sigma(s)dW(s).
\]
Notice that, given a path $X^n_\omega$, 
this equation reduces to equation \eqref{equX} if 
$\chi^n(X_\omega^n)=1$ and
to equation \eqref{equZ} if $\chi^n(X_\omega^n)=0$. 

Now, we construct pathwise the solution process.
For
$\po$-a.e. $\om$,  set  
$ X_\omega^n(t)=X_\omega(t)$ for $0\le t \le \tau^n(\om)$. 
In particular, $\chi^n_{\tau^n(\om)}(X^n_\om)=0$ and
whatever is $X_\omega^n(t)$ for $t>\tau^n(\om)$ 
we will have $\chi^n_t(X_\omega^n)=0$ for $t>\tau^n(\om)$.
Therefore the evolution
of \eqref{eqN} on the time interval $]\tau^n(\om),T]$ 
is given by equation \eqref{equZ}. 
Summing up, 
we have that a solution of equation \eqref{eqN} with the Wiener
process $W$ is the process defined pathwise as follows:
\[X^n(t)=
\begin{cases}
X(t), & t \in [0,\tau^n]\\
\psi^{\tau^n}(X(\tau^n),W)(t), & t \in ]\tau^n,T]
\end{cases}
\]
where  $\psi^{t_0}(y,W)$ 
 denotes the solution of equation \eqref{equZ} (with the Wiener process $W$)
 on the time interval 
$[t_0,T]$ and with initial data $Z(t_0)=y$:
$$ 
 \psi^{t_0}(y,W)(t)=y+\int_{t_0}^t a(s,\psi^{t_0}(y,W))\ ds 
   + \int_{t_0}^t \sigma(s)\ dW(s)
$$
We point out that in this case is enough to assume that 
equation \eqref{equZ} has a unique strong solution
on any time interval $[t_0,T]\subseteq [0,T]$ and that the mapping
$[0,T]\times \mathbb R\times C_0([0,T];\mathbb R ) \ni
(s,y,W) \mapsto \psi^{s}(y,W) \in C([s,T];\mathbb R )$ is
measurable. Actually for fixed initial time, this mapping 
is already known to have nice properties (see, e.g., \cite{IW} Ch. 4 
for the properties of the 
mapping $\psi:\mathbb R\times C_0([0,T];\mathbb R )
\to C([0,T];\mathbb R )$ providing a strong solution on the time
interval $[0,T]$).

An easy example fulfilling these requirements is for the linear 
equation, i.e.  the drift term is $a(t,Z)=c(t)Z(t)$ 
with $c, \sigma$:$[0,T]\to \mathbb R $ measurable and bounded.
Indeed, we have
\begin{equation}\label{caso-lin}
  \psi^{s}(y,W)(t)=e^{\int_s^t c(u)du}y+\int_s^te^{\int^t_u c(r)dr}  \sigma(u)dW(u),
\qquad t\in [s,T].
\end{equation}
\end{remark}

\section{Absolute continuity of $\mu^b$ with respect to $\mu^a$}\label{mains}
We consider equations \eqref{equX} and \eqref{equZ} with the same
initial data $x\in \mathbb R$.
We have the following result. The assumptions are the same as for the 
uniqueness result of the previous section; therefore we denote by
$\mu^b$ the unique law for equation \eqref{equX}.

Let us denote by $\chi_t(Z)$ the indicator function of the set 
$\{\int_0^t \gamma(s,Z)^2 ds<\infty\}$.
\begin{proposition}\label{4.1}
Assume {\bf [A1]} and {\bf [A2]}.\\
If there exists a weak solution
$\big(X, \ (\Omega, \Ff,\{\Ff_t\},\po), \ W\big)$ 
to equation \eqref{equX} such that
\begin{equation}\label{un}\tag{\ref{unic}$^\prime$}
 \po\{\textstyle \int_0^T \gamma(s,X)^2 ds<\infty \}=1,
\end{equation}
 then $\mu^b \prec \mu^a$.
Moreover,
\begin{equation}\label{mu-b-a}
 \frac{d\mu^b}{d\mu^{a}}(Z)=\mathbb E\big[
   e^{\textstyle \mathcal I_T(Z) -\frac 12 \textstyle\int_0^T \gamma(s,Z)^2 ds}
   |\Ff_T(Z)\big]\qquad \po-a.s.,
\end{equation}
where
\begin{equation}\label{rho-lim}
  \mathcal I_T(Z)=\po-\lim_{n \to \infty}  
  \chi_T(Z) \int_0^T \chi^n_s(Z) \gamma(s,Z)dW(s). 
\end{equation}
\end{proposition}
\proof 
Going back to the proof of Proposition \ref{unici},
we have that $\mu^{b,n}\prec \mu^a$ and 
$\|\mu^{b,n}-\mu^b\|_{var}\to 0$.
Then, if $\mu^a(\Lambda)=0$ for some Borelian subset $\Lambda$ of 
$C([0,T];\mathbb R)$, then $\mu^{b,n}(\Lambda)=0$ and finally 
$ \mu^{b}(\Lambda)=0$. This proves $\mu^b\prec \mu^a$.

Moreover,  $\|\mu^{b,n}-\mu^b\|_{var}\to 0$ implies 
that $\mu^{b,n}$ (equivalently, $\po^{*n}$) is a Cauchy
sequence in the metric of total variation.
Since
$\|\po^{*n}-\po^{*m}\|_{var}
=\|\frac{d\po^{*n}}{d\po}-\frac{d\po^{*m}}{d\po}\|_{L^1(\po)}$, 
this is the same as saying that
$$
\frac{d\po^{*n}}{d\po\;\;}= 
 \rho^n_T(Z,W)= e^{\textstyle\int_0^t \chi^n_s(Z) \gamma(s,Z)dW(s) 
-\frac 12 \textstyle\int_0^t \chi^n_s(Z) \gamma(s,Z)^2 ds}
$$ 
is a Cauchy sequence in the metric of $L^1(\po)$. Therefore
$\rho^n_T(Z,W)$ converges in the norm of $L^1(\po)$
to some limit, which is denoted by  $\rho_T(Z,W)$. We want to
identify $\rho_T(Z,W)$.

Notice that if $\int_0^T \gamma(s,Z)^2 ds<\infty$ $\po$-a.s., then 
the stochastic integral in the exponent of $\rho^n_T(Z,W)$
would converge in probability to $\int_0^T  \gamma(s,Z)dW(s)$ 
(see \cite{LS}, Section 4.2.6) and the deterministic integral to 
$\int_0^T \gamma(s,Z)^2 ds$. 
Otherwise, we proceed following the
argument given in \cite{LS} (Section 4.2.9), but with some modification. 
The random variable 
$$ \mathcal I_T^n(Z):=
\chi_T(Z) \int_0^T \chi^n_s(Z) \gamma(s,Z)dW(s)
$$ 
converges in probability. Indeed, 
$ \chi_T(Z)\int_0^T|\chi^n_s(Z) \gamma(s,Z)-\gamma(s,Z)|^2ds$ 
converges to 0 $\po$-a.s.; hence there is convergence in probability.
Therefore, by Lemma 4.6 of \cite{LS},
$\mathcal I_T^n(Z)$ ($n=1,2,\ldots$) is a Cauchy sequence in probability.
It follows that
it converges in probability to a random variabile, which we denote by
$\mathcal I_T(Z)$. 

First, we have
\begin{equation*}
\begin{split}
 \chi_T(Z) \rho_T
&= \chi_T(Z)\
 \po-\lim_n e^{\textstyle\int_0^T \chi^n_s(Z) \gamma(s,Z)dW(s) 
 -\frac 12 \textstyle\int_0^T \chi^n_s(Z) \gamma(s,Z)^2 ds}
\\&= \po-\lim_n \chi_T(Z) 
 e^{\textstyle \int_0^T \chi^n_s(Z) \gamma(s,Z)dW(s) 
 -\frac 12 \textstyle\int_0^T \chi^n_s(Z) \gamma(s,Z)^2 ds}
\\&=
\chi_T(Z) e^{\po-\lim_n[\textstyle\chi_T(Z)  \int_0^T \chi^n_s(Z) \gamma(s,Z)dW(s) 
 -\frac 12 \textstyle\chi_T(Z) \int_0^T \chi^n_s(Z) \gamma(s,Z)^2 ds]}
\\&=
\chi_T(Z) 
e^{\textstyle \mathcal I_T(Z)-\frac 12 \textstyle\int_0^T
  \gamma(s,Z)^2 ds}
\qquad \po-a.s.
\end{split}
\end{equation*}
This means that
\begin{equation}\label{il rho}
 \rho_T(Z,W)=e^{\textstyle \mathcal I_T(Z)-\frac 12 \textstyle\int_0^T
  \gamma(s,Z)^2 ds}  
\end{equation}
a.s. on the set $\{\chi_T(Z)=1\}$.

Now let us check that \eqref{il rho} holds a.s. also on the set 
$\{\chi_T(Z)=0\}$, or equivalently a.s. on the set 
$\{\int_0^T\gamma(s,Z)^2ds=\infty\}$. 
We analyze the left and right hand side of equality
\eqref{il rho}.
The r.h.s. of  \eqref{il rho} vanishes
a.s. on the set $\{\int_0^T \gamma(s,Z)^2 ds=\infty\}$. Indeed, 
by definition $\mathcal I^n_T(Z)=0 $ a.s. on $\{\chi_T(Z)=0\}$ and therefore 
$\mathcal I_T(Z)=0 $ a.s. on $\{\chi_T(Z)=0\}$. 
\\
On the other hand, the l.h.s.  of  \eqref{il rho} vanishes a.s. on the set 
$\{\int_0^T \gamma(s,Z)^2 ds=~\infty\}$. Indeed, on this set
$$
 \rho^n_T(Z,W)=     
e^{\textstyle\int_0^T \chi^n_s(Z) \gamma(s,Z)dW(s) -\frac 12 n}, 
\qquad \po-a.s.
$$
Therefore   
\begin{multline}\label{piccolo}
 \rho^n_T (Z,W)\le e^{-\frac n4} \qquad \po-a.s. \\\;\text{ on } 
\{\textstyle\int_0^T \gamma(s,Z)^2 ds=\infty\}
\cap \{ \textstyle\int_0^T \chi^n_s(Z) \gamma(s,Z)dW(s)\le \frac n4\}.
\end{multline}
Using Chebyshev inequality we get
\begin{equation}\label{cheb}
\begin{split}
\po\{\textstyle\int_0^T \chi^n_s(Z) \gamma(s,Z)dW(s)>\frac n4\} 
&=
\frac 12 \po\{|\textstyle\int_0^T \chi^n_s(Z) \gamma(s,Z)dW(s)|>\frac n4\}
\\
&\le \frac 12 
\frac{\mathbb E[\int_0^T \chi^n_s(Z) \gamma(s,Z)^2 ds]}{(n/4)^2}
\\
&\le \frac 12 \frac n{(n/4)^2}
 \longrightarrow 0 \;\text{ as } n \to \infty.
\end{split}\end{equation}
Let $\chi^{W,n}$ be the indicator function of the set  
$\{\int_0^T \chi^n_s(Z) \gamma(s,Z)dW(s)\le \frac n4\}$.
According to \eqref{cheb} we have that 
\begin{equation}\label{0ins}
 \lim_{n \to \infty} \po\{\chi^{W,n}=0\}=0.
\end{equation}
We investigate the convergence of $\rho^n_T(Z,W)$ on the set
$\{\chi_T(Z)=0\}$:
for any $\varepsilon>0$ we have
\[\begin{split}
\po&\{\rho^n_T(Z,W)[1-\chi_T(Z)]>\varepsilon\}
\\&
=\po\{\rho^n_T(Z,W)>\varepsilon,\ \chi_T(Z)=0\}
\\&
= 
\po\{\rho^n_T(Z,W)>\varepsilon,\ \chi_T(Z)=0,\ \chi^{W,n}=0\}
+
\po\{\rho^n_T(Z,W)>\varepsilon,\ \chi_T(Z)=0,\ \chi^{W,n}=1\}
\\&
\le
\po\{\chi^{W,n}=0\}+
\po\{ e^{-\frac n4}>\varepsilon,\ \chi_T(Z)=0,\ \chi^{W,n}=1\} 
\quad\text{ by } \eqref{piccolo}
\\&
\le 
\po\{\chi^{W,n}=0\}+
\po\{ e^{-\frac n4}>\varepsilon\} 
\longrightarrow 0 \quad\text{  as } n \to \infty \text{ by } \eqref{0ins}.
\end{split}
\]
This implies that $\rho_T(Z,W)[1-\chi_T(Z)]=0$ a.s.;
hence,  $\rho_T(Z,W)=0$ a.s. on the set $\{\chi_T(Z)=0\}$.

We conclude that
\begin{equation}\label{insiemi}
 \rho_T(Z,W)=e^{\textstyle \mathcal I_T(Z)-\frac 12 \textstyle\int_0^T
  \gamma(s,Z)^2 ds} \qquad \po-\text{ a.s. }
\end{equation}

Finally, denoting by $\po^*$ the limit of $\po^{*n}$, so that
$\mu^b(\Lambda)=\po^*\{Z\in \Lambda\}$, we have proved that
$\frac{d\po^*}{d\po\;}=\rho_T(Z,W)$.
As done in the proof of Theorem \ref{semplice}, 
we get \eqref{mu-b-a}.
\hfill$\Box$

\section{Equivalence of the laws}\label{s-equi}
As noticed in the previous section, if
$\po\{\textstyle \int_0^T \gamma(s,Z)^2 ds<\infty \}=1$,
then 
\begin{equation}\label{casoINT}
 \mathcal I_T(Z)=\int_0^T \gamma(s,Z)dW(s)
\end{equation}
and therefore $\dfrac{d\po^*}{d\po\;}=\rho_T(Z,W) $, where 
$$
 \rho_t(Z,W) 
 =e^{\textstyle \int_0^t \gamma(s,Z)dW(s) -\frac 12 \textstyle\int_0^t \gamma(s,Z)^2 ds}
$$
is a strictly positive martingale.

From this, we have a result on how to use Girsanov transform under
very weak assumptions (basically, avoiding Novikov condition or
similar conditions involving the expectation of the exponential of a
random variable related to the integral of $\gamma(s,Z)$;
see \cite{nov}, \cite{Ka}, \cite{K}).
\begin{theorem}\label{p1}
Assume {\bf [A1]}, {\bf [A2]} and that  equation \eqref{equX}
has a weak solution 
$ \big(X, \ (\Omega, \Ff,\{\Ff_t\},\po), \ W\big) $.
Denote by $Z$ the unique solution of equation \eqref{equZ} 
with respect to the same stochastic basis and  Wiener process.\\
If
\begin{equation}\label{primo-u}
 \po\{\textstyle \int_0^T \gamma(s,X)^2 ds<\infty \}=1,
\end{equation}
\begin{equation}\label{primo-z}
 \po\{\textstyle \int_0^T \gamma(s,Z)^2 ds<\infty \}=1,
\end{equation}
then\\
i) the process $\rho=\rho(Z,W)$ given by
\begin{equation}\label{dens}
\rho_t=
e^{\textstyle\int_0^t\gamma(s,Z)dW(s)
 -\frac 12 \textstyle\int_0^t \gamma(s,Z)^2 ds}, \quad
0\le t \le T, 
\end{equation}
is a positive  $\{\Ff_t\}$-martingale; 
in particular
\begin{equation}\label{mart}
 \mathbb E [\rho_t(Z,W)] = 1  \;\text{ for any } t \in [0,T].
\end{equation}
ii)
\begin{equation}\label{w.star}
W^*(t)=W(t)-\int_0^t  \gamma(s,Z)\ ds \ , \qquad t \in [0,T],
\end{equation}
is a Wiener process with respect to
$\po^*$,
where
the probability measure $\po^*$ is defined on $(\Omega,\Ff_T)$ by 
\begin{equation}\label{lap*}
 d\po^* = \rho_T(Z,W) \  d\po.
\end{equation}
\end{theorem}
\proof 
$i)$
Notice that the exponential process
 $\rho(Z,W)$ is a positive local martingale and then, by Fatou lemma, 
a supermartingale. Since $\rho_0(Z,W)=1$, it is enough
to have $\mathbb E [\rho_T(Z,W)]=1$ 
 in order to prove that it is a martingale. But, $\rho_T(Z,W)$  is the
$L^1(\po)$-limit of $\rho^n_T(Z,W)$; since  
we already know from the proof of Proposition \ref{unici} that
$$
\mathbb E [\rho^n_T(Z,W)]=1 \quad \text{
  for any } n=1,2,\dots 
$$
we get that $\mathbb E [\rho_T(Z,W)]=1$.

$ii)$  Given $i)$, this is Girsanov theorem
(see, e.g., \cite{gir}).
\hfill$\Box$

\smallskip

Now we state our main result.
\begin{theorem}\label{main}
Assume {\bf [A1]}, {\bf [A2]} and that  equation \eqref{equX}
has a weak solution \\
$ \big(X, \ (\Omega, \Ff,\{\Ff_t\},\po), \ W\big) $.
Denote by $Z$ the unique solution of equation \eqref{equZ} 
with respect to the same stochastic basis and  Wiener process.\\
If \eqref{primo-u}-\eqref{primo-z} hold,
then the law of the solution of equation \eqref{equX} is unique.
Moreover, $\mu^b \sim \mu^a$. In particular, 
\begin{alignat}{1}
\frac{d\mu^b}{d\mu^a}(Z)=
  \mathbb E
  \left[ e^{+\int_0^T \gamma(s,Z)dW(s)-\frac 12 \int_0^T \gamma(s,Z)^2 ds}
   \big|\Ff_T(Z)\right]  \qquad\quad 
 &\po-a.s. \label{roz}  \\
\frac{d\mu^a}{d\mu^b}(Z)=
 \mathbb E
 \left[e^{-\int_0^T \gamma(s,Z)dW(s)+\frac 12 \int_0^T \gamma(s,Z)^2 ds} |\Ff_T(Z)\right]
 \qquad\quad&\po-a.s.  \label{leggi5} \\
\frac{d\mu^a}{d\mu^b}(Z)=
 \mathbb E^*
 \left[e^{-\int_0^T \gamma(s,Z)dW^*(s)
       -\frac 12 \int_0^T \gamma(s,Z)^2 ds} |\Ff_T(Z)\right]
 \quad\quad&\po^*-a.s. \label{densitaleggi2}\\
\frac{d\mu^a}{d\mu^b}(X)=
  \mathbb E
  \left[ e^{-\int_0^T \gamma(s,X)dW(s)-\frac 12 \int_0^T \gamma(s,X)^2 ds}
   \big|\Ff_T(X)\right]  \quad\quad 
  &\po-a.s. \label{densX} 
\end{alignat}
where $\po^*, W^*$ are defined by \eqref{lap*}, \eqref{w.star}
respectively.
\end{theorem}
\proof 
Uniqueness in law comes from Proposition \ref{unici}, $\mu_b\prec \mu_a$
from Proposition \ref{4.1} and \eqref{roz}
from \eqref{mu-b-a}, \eqref{casoINT} with the assumption \eqref{primo-z}.

Moreover,  \eqref{primo-z} implies that 
$\po\{\rho_T(Z,W)=0\}=0$. Then $\po\prec \po^*$
with 
$$
 \frac{d\po\;}{d\po^*}= \frac 1{\rho_T(Z,W)}, \qquad \po^*-a.s.
$$
(see Lemma 6.8 in \cite{LS}).
From $\po \prec \po^*$ follows $\mu^a \prec \mu^b$. 

As done in the proof of Theorem \ref{semplice}, 
from $\frac{d\po\;}{d\po^*}=\big(\rho_T(Z,W)\big)^{-1}$ we get \eqref{leggi5}.
Moreover, using \eqref{w.star} we get
$$
 \frac{d\po\;}{d\po^*}=
 e^{-\int_0^T \gamma(s,Z)dW(s)+\frac 12 \int_0^T \gamma(s,Z)^2 ds}
 =e^{-\int_0^T \gamma(s,Z)dW^*(s)-\frac 12 \int_0^T \gamma(s,Z)^2 ds};
$$
in the same way, this gives \eqref{densitaleggi2}.
In particular, 
\begin{equation}
 \mathbb E^*
 \Big[e^{\textstyle-\int_0^T \gamma(s,Z)dW^*(s)
       -\frac 12 \int_0^T \gamma(s,Z)^2 ds} \Big]=1.
\end{equation}
This is written for the solution 
$ \big(Z, \ (\Omega, \Ff,\{\Ff_t\},\po^*), \ W^*\big)$ 
of equation \eqref{equX}. Since there is uniqueness in law for
equation \eqref{equX}, if we consider the same relationship for
the solution $ \big(X, \ (\Omega, \Ff,\{\Ff_t\},\po), \ W\big)$
we get 
\begin{equation}\label{1X}
 \mathbb E \Big[e^{\textstyle-\int_0^T \gamma(s,X)dW(s)
       -\frac 12 \int_0^T \gamma(s,X)^2 ds} \Big]=1.
\end{equation}
Now, let us start from equation \eqref{equX} and 
consider equation \eqref{equZ} as a modification of equation
\eqref{equX} by a change of drift.
Thanks to \eqref{1X}, we can use Girsanov theorem.
Similarly to what we have done in Theorem \ref{semplice} and Remark
\ref{ledensi} ii),
we get \eqref{densX}. 
\hfill$\Box$

\section{Conclusions} \label{conclu}
We compare our results with whose of Liptser and Shiryaev.
Let us begin reminding a result of \cite{LS}: there Theorem 7.18
states that if we assume {\bf [A1]}, {\bf [A2]}, 
that equation \eqref{equX} has a
weak solution $\big(X,\ (\Omega, \Ff,\{\Ff_t\},\po), \ W \big)$
satisfying \eqref{primo-u} and
\begin{equation}\label{novk}
\mathbb E\Big[e^{\textstyle -\int_0^T \gamma(s,X)dW(s)
   -\frac 12 \int_0^T \gamma(s,X)^2 ds} \Big]=1,
\end{equation}
then $\mu^b \sim \mu^a$ and 
$$
 \frac{d\mu^a}{d\mu^b}(X)=
  \mathbb E\Big[  e^{\textstyle -\int_0^T \gamma(s,X)dW(s)
   -\frac 12 \int_0^T \gamma(s,X)^2 ds}
  \big|\Ff_T(X)\Big], \qquad
 \po-a.s.
$$
The crucial issue is how to get \eqref{novk} without assuming the 
quite strong Novikov condition (see \cite{nov}) 
$$
 \mathbb E\Big[e^{\textstyle \frac 12 \int_0^T \gamma(s,X)^2 ds} \Big]<\infty
$$
or other conditions 
involving the expectation of the exponential of a
random variable related to the integral of $\gamma(s,X)$
(see \cite{Ka}, \cite{K}). This is done in our Theorem \ref{main} with
the ''$\po$-a.s.'' conditions \eqref{primo-u}-\eqref{primo-z}.

However, Liptser and Shiryaev present another result, more operative
than Theorem 7.18. This is 
Theorem 7.19 of \cite{LS} providing $\mu^b \sim \mu^a$
with the same assumptions of Theorem 7.18 except \eqref{novk},
which is replaced by
\begin{equation}\label{loro}
 \po\{\textstyle \int_0^T \gamma_a(s,X)^2  ds <\infty \}=
       \po\{\textstyle \int_0^T \gamma_b(s,X)^2 ds <\infty \}=1,
\end{equation}
\begin{equation}\label{loro2}
\po\{\textstyle \int_0^T \gamma_a(s,Z)^2  ds <\infty \}=
   \po\{\textstyle \int_0^T \gamma_b(s,Z)^2 ds <\infty \}=1,
\end{equation}
where 
$$ 
 \gamma_a(s,X)= \sigma^+(s,X) a(s,X), \qquad
 \gamma_b(s,X)=\sigma^+(s,X) b(s,X).
$$
Because of $\gamma=\gamma_b-\gamma_a$, assumptions
\eqref{loro}-\eqref{loro2} are stronger than \eqref{primo-u}-\eqref{primo-z}.

Therefore, we can see our Theorem \ref{main} as an intermediate result  
between the two theorems of \cite{LS}. We have the same result as
Theorem 7.18, but saying concretely how  to get \eqref{novk} with
''$\po$-a.s.'' conditions. This is
in the same spirit as Theorem 7.19. However, our
conditions \eqref{primo-u}-\eqref{primo-z} 
involve only the difference $b-a$ of the  drift terms, whereas 
conditions \eqref{loro}-\eqref{loro2} involve  both the
drift terms $b$ and $a$.

We point out that our results on the absolute continuity 
of the laws are identical to \cite{LS}, but the expressions of 
the Radon-Nikodym derivatives are different
from those of  Liptser and Shiryaev.
In fact, under \eqref{loro}-\eqref{loro2}
Theorem 7.19 gives
\begin{equation}\label{noW}
 \frac{d\mu^b}{d\mu^a}(Z)=
 e^{\int_0^T \sigma^+(s,Z)^2[b(s,Z)-a(s,Z)]\ dZ(s)
-\frac 12 \int_0^T \sigma^+(s,Z)^2[b(s,Z)^2-a(s,Z)^2]\ ds} 
\end{equation}
and
\begin{equation} \label{noW2}
 \frac{d\mu^a}{d\mu^b}(X)=
 e^{-\int_0^T \sigma^+(s,X)^2[b(s,X)-a(s,X)] dX(s)
+\frac 12 \int_0^T \sigma^+(s,X)^2[b(s,X)^2-a(s,X)^2] ds} 
\end{equation}
($\po$-a.s.). Let us show that \eqref{noW} can be obtained from
\eqref{roz};
similarly, for \eqref{noW2} from \eqref{densX}.
If this is true, then we conclude that 
with our proofs we get the same result as Liptser and Shiryaev.
However, our proofs  are different from  \cite{LS}
basically in one point: Liptser and Shiryaev analyze the equation
satisfied by the Radon-Nikodym derivative,  whereas we analyze the 
Radon-Nikodym derivative as the limit of the sequence $\frac{d\po^{*n}}{d\po}$.
This makes our proofs shorter.

Let us come back to the expression of the Radon-Nikodym derivative 
$\frac{d\mu^b}{d\mu^a}(Z)$. Now, we assume
\eqref{loro}-\eqref{loro2};
of course our results hold true. Therefore $\frac{d\mu^b}{d\mu^a}(Z)$
is given by \eqref{roz}.
From \eqref{equZ} we have $W$ depending on $Z$: 
$dW(t)=\sigma^+(t,Z)[dZ(t)-  a(t,Z)\ dt]$.
Then  
$$
\int_0^T \gamma(s,Z)dW(s)-\frac 12 \int_0^T \gamma(s,Z)^2 ds 
$$
becomes (formally)
\begin{equation}\label{LSden}
 \int_0^T \sigma^+(s,Z)^2[b(s,Z)-a(s,Z)] dZ(s)
  -\frac 12 \int_0^T \sigma^+(s,Z)^2[b(s,Z)^2-a(s,Z)^2] ds.
\end{equation}
Since
\[
 \int_0^T |\sigma^+(s,Z)^2[b(s,Z)-a(s,Z)]|^2 \sigma(s,Z)^2 \ ds 
\le \int_0^T \gamma(s,Z)^2 \ ds ,
\]
the stochastic integral  in \eqref{LSden} 
is well defined if \eqref{primo-z} holds, whereas
the deterministic integral requires \eqref{loro2}.
Then, \eqref{LSden} is  in fact well
defined with assumption \eqref{loro2} and it  
depends only on $Z$ and not also on $W$. 
From \eqref{roz} and \eqref{LSden}, we get \eqref{noW}.

\section{Bigger dimensions}\label{bigd}
Let $d, m \in \mathbb N$ with $dm>1$. 
The solution processes $X$ and $Z$ 
have paths in $C([0,T];\mathbb R^d)$, 
the 
initial data $x\in \mathbb R^d$ and $W$ is an $m$-dimensional Wiener
process. Any  vector $v$ is a  column vector, whose transposed is the
row vector $v^T$. We set
$\|X\|^2=\sum_{i=1}^dX_i^2$.

Let $\mathcal B_t$ be the $\sigma$-algebra of Borelian subsets of
$C([0,t];\mathbb R^d )$, for $0<t\le T$.
The drift terms  $a$ and $b$ are $\mathbb R^d$-valued non 
anticipative measurable functionals, that is
$$
 a, b: [0,T]\times C([0,T];\mathbb R^d)\to \mathbb R^d, 
$$
are measurable and, for each $t \in [0,T]$,  $a(t,\cdot), b(t,\cdot)$ 
are  $\mathcal B_{t}$-measurable.
Similarly, 
the diffusion term  
$\sigma:[0,T]\times C([0,T];\mathbb R^d)\to \mathbb R^d\times \mathbb R^m$
is a non anticipative measurable functional;
in particular
$(\sigma W)_i=\sum_{k=1}^m\sigma_{ik} W_k$ for $i=1,\ldots, d$.
The entries satisfy
the previous conditions;
the two main assumptions  become
\[
\hspace{-4mm}{\mathbf{ [A1]}}
\left[
\begin{array}{l}
\; \exists \text{ constants } L_1,L_2 \text{ and a function }
K \text{ non decreasing and right continuous,}\\
\text{ with }0 \le K(s)\le 1, 
\text{ such that all the components  }
a_i, b_{ik} \text{ satisfy }
\\[2mm]
 a_i(t,Y)^2+\sigma_{ik}(t,Y)^2 
\le L_1 \int_0^t [1+\|Y(s)\|^2] dK(s)+L_2[1+\|Y(t)\|^2]
  \\\hspace*{66mm} \forall t \in [0,T], Y \in C([0,T];\mathbb R^d)
\\\text{ and } \\
 |a_i(t,Y_1)-a_i(t,Y_2)|^2+|\sigma_{ik}(t,Y_1)-\sigma_{ik}(t,Y_2)|^2
\\
\hspace*{30mm}
\le L_1\int_0^t \|Y_1(s)-Y_2(s)\|^2dK(s)+L_2 \|Y_1(t)-Y_2(t)\|^2 
 \\\hspace*{66mm} \forall t \in [0,T], Y_1,Y_2 \in C([0,T];\mathbb R^d)
\end{array}
\right.
\]
\[
\hspace{-8mm}{\mathbf{ [A2]}}
\left[
\begin{array}{l}
\; \exists \gamma \text{ finite and } \mathbb R^m\text{-valued 
non anticipative measurable functional: }\\ 
\quad\sigma(s,Y)\gamma(s,Y)=b(s,Y)-a(s,Y)\qquad
\forall s\in [0,T], Y \in C([0,T];\mathbb R^d).
\end{array}
\right.
\]
In Remark \ref{sol:deb}, the solution mapping 
is 
$$
[0,T]\times \mathbb R^d\times C_0([0,T];\mathbb R^m ) \ni
(s,y,W) \mapsto \psi^{s}(y,W) \in C([s,T];\mathbb R^d )
$$
and the linear equation in the example has solution
still given by \eqref{caso-lin},
where $c:[0,T]\to\mathbb R^d\times \mathbb R^d$ is measurable and bounded.

All the results of the previous sections hold true with the suitable
change of notations. Mainly, $\gamma^2$ becomes 
$\|\gamma\|^2= \sum_{k=1}^m\gamma_k^2$ and
$\sigma^+$ is the pseudo-inverse matrix of $\sigma$ (see, e.g.,
\cite{matrici}, \cite{LS2}); $\sigma^+$ is an $m\times d$-matrix, 
uniquely defined.

However, let us investigate this multidimensional problem in details.
Assumption {\bf [A2]} refers to the linear system of 
$d$ equations in $m$ unknowns
\begin{equation}\label{sg}
\sigma(s,Y)\gamma(s,Y)=b(s,Y)-a(s,Y)
\end{equation} 
and is a consistency condition involving $\sigma$ and $b-a$ 
(see, e.g., \cite{matrici} for all the results on linear systems and 
matrices). Moreover, 
if a solution $\gamma$ exists and
 $\operatorname{rank} \sigma=m$  
then the solution  of \eqref{sg} is unique and is given by 
\begin{equation} \label{Gamma} 
\gamma(s,Y)=\sigma^+(s,Y)[b(s,Y)-a(s,Y)].
\end{equation}
In particular, if $\sigma$ has maximal rank we have
\begin{equation}\label{s+}
 \sigma^+= \sigma^T(\sigma\sigma^T)^{-1} 
   \qquad \text{ or }\;
 \sigma^+= (\sigma^T \sigma)^{-1} \sigma^T.
\end{equation}
(Notice that if the rank of $\sigma$ is maximal, then 
also the square matrices $\sigma \sigma^T$ and $\sigma^T\sigma$
have maximal rank and therefore are invertible.)
And if the matrix $\sigma$ vanishes, then also $\sigma^+$ has all the
entries equal to 0.
\\
Otherwise, there are infinite solutions of \eqref{sg}, 
one of them being given by 
\eqref{Gamma}. This is the case of dimensions $d \ge m$ with rank of 
$\sigma$ not maximal ($ <m$), or of dimensions
$m>d$. To handle these cases,
let us
recall the singular value decomposition of the 
$d\times m$-matrix $\sigma$
with $\operatorname{rank} \sigma=r$, $r\le \min(d,m)$:
$$
 \sigma=\lambda^{(1)}  u^{(1)} (v^{(1)})^T+\lambda^{(2)}   u^{(2)} (v^{(2)})^T
   +\ldots  +\lambda^{(r)}  u^{(r)} (v^{(r)})^T,
$$
where $\lambda^{(1)}\ge \lambda^{(2)}\ge \ldots \lambda^{(r)}>0$,
$\{u^{(i)}\}_{i=1}^r$ is an orthonormal set of $d$-dimensional vectors 
and $\{v^{(i)}\}_{i=1}^r$  is an orthonormal set of $m$-dimensional
vectors. Moreover
$$
 \sigma^+=\frac 1{\lambda^{(1)}} v^{(1)} (u^{(1)})^T
         +\frac 1{\lambda^{(2)}} v^{(2)} (u^{(2)})^T
   +\ldots  +\frac 1{\lambda^{(r)}}  v^{(r)} (u^{(r)})^T.
$$
Then the $d$-dimensional vector $\sigma(t,X) dW(t)$
can be written as
$$
 \lambda^{(1)}(t,X) u^{(1)}(t,X) v^{(1)}(t,X)^TdW(t) 
 +\ldots    
 +\lambda^{(r)}(t,X) u^{(r)}(t,X) v^{(r)}(t,X)^T dW(t).
$$
This means that  
$$
\sigma(t,X) dW(t)=  
\tilde \sigma(t,X) d\tilde W^X(t)  
$$
with
$\tilde \sigma(t,X)$ the $d\times r$-matrix  and $\tilde W^X(t)$ the $r$-vector
defined by
\[\begin{split}
 &\tilde \sigma_{ij}(t,X)=\lambda^{(j)}(t,X) u^{(j)}_i(t,X), \\ 
 & \tilde W^X_i(t)= \int_0^t v^{(i)}(s,X)^T dW(s)
  \equiv \sum_{k=1}^m \int_0^t v_k^{(i)}(s,X) dW_k(s).
\end{split}\]
The matrix $\tilde \sigma(t,X)$ has maximal rank ($= r$), since 
the vectors $u^{(j)}$ are orthogonal to each other. 
Moreover, 
$\tilde W^X$ is an $r$-dimensional Wiener
process. Indeed, the components of this vectors are one dimensional
independent Wiener processes, thanks to the fact that
$\{v^{(i)}(s,X)\}_{i=1}^r$  is an orthonormal set of $m$-dimensional
vectors.

Now, let us read  equation
\eqref{equX} 
with $\tilde \sigma(t,X) d\tilde W^X(t)$ 
instead of $\sigma(t,X) dW(t)$; similarly for equation \eqref{equZ}. 
Then, according to the previous considerations
for the $d\times r$-diffusion matrix $\tilde \sigma$
with $r<d$ and   maximal rank ($=r$),  
we get that the system
$$
 \tilde \sigma(s,Y)\tilde \gamma(s,Y)=  b(s,Y)- a(s,Y) ;  
$$
has at most one solution, given by
$$
 \tilde \gamma(s,Y)=\tilde \sigma^{+}(s,Y)[b(s,Y)- a(s,Y)]  
$$
with
$ \tilde \sigma^{+}=(\tilde \sigma^T \tilde \sigma)^{-1} \tilde \sigma^T$.
Actually, $\tilde \sigma$ has $r$ columns given by the vectors
$\lambda^{(i)} u^{(i)}$ and 
$\tilde \sigma^+$ has $r$ rows given by the vectors 
$(\frac 1{\lambda^{(i)}} u^{(i)})^T$.

A tedious but easy computation provides
$$
 \|\sigma^{+}(b-a)\|^2 
 = \sum_{i=1}^r\left|\frac {{u^{(i)}}^T(b-a)}{\lambda^{(i)}}\right|^2
 =\| \tilde  \sigma^+(b-a)\|^2.
$$
Since the latter
quantity is uniquely defined, also the first is unique.
Therefore 
$$\int_0^T\|\gamma(s,Y)\|^2 ds=\int_0^T\|\tilde
\gamma(s,Y)\|^2ds \quad \text{ for } Y \in C([0,T];\mathbb R^d),
$$ 
where
$$
 \gamma(s,Y)= \sigma^{+}(s,Y)[b(s,Y)- a(s,Y)]. 
$$
This expression of $\gamma$ provides 
the  unique relevant solution of \eqref{sg}
 in the Girsanov transform, even when the solution of \eqref{sg} 
is not unique. In particular, we have
$$
\frac{d\mu^b}{d\mu^a}(Z)=  
 \mathbb E
 \Big[ e^{\textstyle +\int_0^T \gamma(s,Z)dW(s)-\frac 12 \int_0^T \|\gamma(s,Z)\|^2 ds}
  \Big|\Ff_T(Z)\Big] , 
$$
$$
\frac{d\mu^a}{d\mu^b}(X)=  
 \mathbb E 
 \Big[e^{\textstyle 
 -\int_0^T\gamma(s,X)dW(s)-\frac 12 \int_0^T \|\gamma(s,X)\|^2 ds}
 \Big|\Ff_T(X)\Big],
$$ 
a.s., when we assume
$$
\po\{\textstyle \int_0^T \|\gamma(s,X)\|^2 ds<\infty \}=1,\quad
 \po\{\textstyle \int_0^T \|\gamma(s,Z)\|^2 ds<\infty \}=1,
$$
instead of \eqref{primo-u}-\eqref{primo-z}.

We therefore conclude that we get all our previous results, included the
uniqueness result. Let us emphasize  that the
uniqueness question is not investigated in \cite{LS} 
for $dm>1$ when equation \eqref{sg} has more than
one solution.
However, even if not stated, it appears clear from the results of
\cite{LS}  in the one
dimensional case  that there is uniqueness in law for equation
\eqref{equX},
because of the uniqueness of $\gamma$ (see also the beginning of
Section \ref{coro}).

\section{Applications}\label{sAPP}
Let us consider the case of 
 $b=a+f$, that is we deal
with
\begin{align*}
&dX(t)=a(t,X)\  dt+f(t,X)\  dt+\sigma(t,X)\ dW(t), 
 \quad  &X(0)=x\\
&dZ(t)=a(t,Z)\ dt+\sigma(t,Z)\ dW(t), \quad  &Z(0)=x
\end{align*}
To apply the results of \cite{LS}, besides {\bf [A1]} and {\bf [A2]}
we have to check
conditions on $a$ and $a+f$, whereas our results require only a
condition on $f$.  
Let us see how to use our results; first, in the one
dimensional problem, then in the infinite dimensional one.
\\[3mm]
{\bf One dimensional stochastic differential equations}
\nopagebreak

We consider conditions involving the process $X$; of course, 
the same holds true
for those involving $Z$.

Our condition \eqref{primo-u} becomes
\begin{equation}\label{pprr}
  \po\{\textstyle \int_0^T \sigma^+(t,X)^2f(t,X)^2\ dt <\infty \}=1,
\end{equation}
whereas \eqref{loro} becomes
\begin{equation}\label{somma2}
 \po\{\textstyle \int_0^T \sigma^+(t,X)^2
 \big(a(t,X)^2+[a(t,X)+f(t,X)]^2\big)\ dt <\infty \}=1,
\end{equation}
that is 
\begin{multline*}
\po\{\textstyle \int_0^T \sigma^+(t,X)^2 a(t,X)^2 \ dt<\infty\}
\\=1=
\po\{\textstyle \int_0^T \sigma^+(t,X)^2[2a(t,X)f(t,X)+f(t,X)^2]\ dt
<\infty \}.
\end{multline*}
If $af\ge 0$ this is equivalent to 
\begin{equation}\label{rrpp}
 \po\{\textstyle \int_0^T \sigma^+(t,X)^2a(t,X)^2\ dt <\infty\}=
 \po\{\textstyle \int_0^T \sigma^+(t,X)^2f(t,X)^2\ dt <\infty \}=1,
\end{equation}
In general, the latter implies \eqref{somma2}.

This condition \eqref{rrpp}  is stronger than \eqref{pprr},
{\it unless}  $\sigma$ is constant. In fact, if $\sigma$ is a constant
$\neq 0$, 
then $\po\{\int_0^T \sigma^+(t,X)^2 a(t,X)^2 \ dt<\infty\}=1$
becomes
\begin{equation}\label{solo-a}
 \po\{\textstyle \int_0^T a(t,X)^2 \ dt < \infty \}=1
\end{equation}
which is trivially fulfilled thanks to the growth condition on $a$ included
in {\bf [A1]}. Then we only have to check if
$$\po\{\textstyle \int_0^T f(t,X)^2\ dt <\infty \}=1,$$
that is \eqref{pprr} and \eqref{rrpp} are equivalent.

Otherwise, for general $\sigma$, condition \eqref{rrpp} is stronger than 
our condition \eqref{pprr}. 
\\[3mm]
{\bf Infinite dimensional stochastic differential equations}

To have  weaker assumption is even more important 
in the infinite dimensional setting. In fact,
{\it even if} $\sigma$ is constant, 
the conditions of Liptser and Shiryaev 
\begin{equation}\label{nonv}
\begin{split}
 \po\{\textstyle \int_0^T \|\sigma^+(s,X)a(s,X)\|^2\ ds <\infty\}=1\\
 \po\{\textstyle \int_0^T \|\sigma^+(s,Z)a(s,Z)\|^2\ ds <\infty\}=1
\end{split}
\end{equation}
may be cumbersome (see next Remark \ref{oss}). This is different from
the finite dimensional framework. Indeed, the coefficients
 $\sigma$ and $a$ are now operators in some infinite dimensional
spaces.
Our results allow to obtain uniqueness in law and absolute continuity
of the laws getting rid of \eqref{nonv}.

First, we fix the Hilbert spaces to work in and we make precise the
norm to consider in \eqref{nonv}.
We are given  separable Hilbert spaces 
$E\subseteq E_1\subseteq H$ with continuous and dense embeddings.
The space $E$ will ''replace'' the state space $\mathbb R^d$.
We denote by $\|\cdot \|_H$ the norm in $H$
and
by $\!_H\langle\cdot,\cdot\rangle_H$ the scalar
product in $H$.

For simplicity, let us consider the very simple but interesting case
of constant diffusion and drift independent of the first variable $t$
and linear in the second variable $X$. 
Equation \eqref{equZ} becomes
\begin{equation}\label{ou}
 dZ(t)=AZ(t)\ dt+\sqrt Q \ dW(t), \qquad Z(0)=x
\end{equation}
where $W$ is 
 a cylindrical Wiener process in $H$,
defined on the probability space 
$(\Omega, \Ff,\{\Ff_t\},\po)$.
This means that, if $\{e_j\}_{j=1}^\infty$ is a complete
orthonormal system of $H$, then we represent $W(t)=\sum_j
\beta_j(t) e_j$ with $\{\beta_j\}_{j=1}^\infty$  
a sequence of  i.i.d. one dimensional Wiener processes defined on
$(\Omega, \Ff,\{\Ff_t\},\po)$, .
The operators $A$ and $Q$ are linear operators in $H$ and $x \in E$.
Therefore equation \eqref{ou} is a linear stochastic equation; this is
the simplest infinite dimensional equation to deal with, for which 
it is easy to get existence and uniqueness of solutions and of invariant
measures (see, e.g., \cite{dpz}, \cite{dpz2}). More general equations
can be deal with in a similar way; but \eqref{ou} allows us already to cover
interesting examples.

Instead of \eqref{equX}, consider the semilinear stochastic equation
\begin{equation}\label{seml}
 dX(t)=[AX(t)+F(X(t))]\ dt+\sqrt Q dW(t), \qquad X(0)=x
\end{equation}
where $F:E\to E_1$ is measurable.

According to Remark \ref{sol:deb}, we assume 
\[
{\mathbf{ [A3]}}
\left[
\begin{array}{l}
\text{for 
 any initial data $x \in E$ and on any  time interval }
 [t_0,T]\subseteq [0,T]\\
\text{equation \eqref{ou}
has a unique strong solution } Z, 
 \text{ whose paths are in }\\
C([t_0,T];E) \text{ a.s.} 
\end{array}
\right.
\]
General conditions on $A$ and $Q$ to get it, can be found, e.g., in
\cite{dpz}, whereas examples are in \cite{dpz} and \cite{cosa}.
\\
Moreover \medskip

$
{\mathbf{ [A4]}}
\left[
\begin{array}{l}
Ran(F)\subseteq Ran(\sqrt Q) \\ \exists\  (\sqrt Q)^{-1}
\end{array}
\right.
$
\\
We set
$$
 \Gamma(Y)=(\sqrt Q)^{-1} F(Y) \qquad \forall Y\in E.
$$
We have that $\Gamma:E\to H$   is measurable.

Here is our  result of uniqueness in law of Section \ref{coro}, 
stated in the infinite dimensional setting.
\begin{proposition}
Assume {\bf [A3]} and {\bf [A4]}.\\
If there exist two weak solutions
$\big(X, (\Omega, \Ff,\{\Ff_t\},\po), \ W\big)$ and 
$\big(X^\prime, (\Omega^\prime, \Ff^\prime,\{\Ff^\prime_t\},\po^\prime), 
\ W^\prime\big)$ to equation \eqref{seml}  
 with the same initial data $x \in E$ and  with paths in 
$C([0,T];E)$ $\po$-a.s., 
such that
\begin{equation}
 \po\{\textstyle \int_0^T \|\Gamma(X(s))\|_H^2 ds<\infty
 \}=
 \po^\prime\{\textstyle \int_0^T \|\Gamma(X^\prime(s))\|_H^2 ds<\infty
 \}=1,
\end{equation}
 then the laws of $X$ and $X^\prime$ are the same.
\end{proposition}

For the equivalence of the laws we have
\begin{theorem}\label{infEquiv}
Assume  {\bf [A3]} and {\bf [A4]}. \\
Given $x \in E$, if there exists  a weak solution 
$\big(X, (\Omega, \Ff,\{\Ff_t\},\po), \ W\big)$
to equation \eqref{seml} with paths in 
$C([0,T];E)$ $\po$-a.s., and satisfying 
\[
 \po\{\textstyle \int_0^T\|\Gamma(X(s))\|_H^2 ds<\infty \}=1,
\]
then there is uniqueness in law for equation \eqref{seml} and 
$\mu^{A+F}\prec \mu^A$. Further, if the strong solution
$\big(Z, (\Omega, \Ff,\{\Ff_t\},\po), \ W\big)$ to equation \eqref{ou}
satisfies
\[
  \po\{\textstyle \int_0^T\|\Gamma(Z(s))\|_H^2 ds<\infty \}=1,
\]
then $\mu^{A+F}\sim \mu^A$; in particular
\[
 \frac{d\mu^{A+F}}{d\mu^A\;}(Z)=
\mathbb E\Big[ e^{\textstyle 
   + \int_0^T \!\!\!_H\langle  \Gamma(Z(s)),dW(s)\rangle_H
   -\frac 12 \int_0^T
  \|\Gamma(Z(s)\|_H^2ds} \big|\Ff_T(Z)\Big], 
\]
\[
 \frac{d\mu^{A}\;}{d\mu^{A+F}}(X)=
\mathbb E\Big[ e^{\textstyle 
   - \int_0^T \!\!\!_H\langle  \Gamma(X(s)),dW(s)\rangle_H
   -\frac 12 \int_0^T
  \|\Gamma(X(s)\|_H^2ds} \big|\Ff_T(X)\Big],
\]
$\po$-a.s.
\end{theorem}

\begin{remark}\label{oss}
Conditions \eqref{nonv} become
\[
\po\{\textstyle\int_0^T \|(\sqrt Q )^{-1}AX(s)\|_H^2 \ ds<\infty\}=1
\]
\[
\po\{\textstyle\int_0^T \|(\sqrt Q )^{-1}AZ(s)\|_H^2 \ ds<\infty\}=1.
\]
We point out that they are
not satisfied  in  the example of Section 4 in \cite{cosa}.
\end{remark}

\medskip 
We conclude analyzing a consequence  of the equivalence of the laws.
If $\mu^{A+F}\sim \mu^A$, then for any fixed $t\in [0,T]$
  the law of $X(t)$ 
is equivalent to the law of $Z(t)$.
If we know properties of the law of $Z(t)$, then they hold a.s. also
for $X(t)$. 
Usually,  properties of the solutions of the linear equation \eqref{ou}
are easier to obtain than for the non linear equation \eqref{seml}.
The Girsanov transform allows to link these results.
An important application is in the  study of the
asymptotic behaviour, as $t \to \infty$,  of equation \eqref{seml}
in an infinite dimensional space (see, e.g., \cite{dpz2} for the  general
theory and examples, and  \cite{cosa} for examples) when our
results hold on any finite time interval $[0,T]$, that is for any $T>0$.

\bibliographystyle{amsplain}

\begin{thebibliography}{11}

\bibitem{matrici}
Campbell, S. L., Meyer, C. D.:
{\it Generalized inverses of linear transformations},
Dover, 1991. 

\bibitem{che}
Cherny, A. S.:
On the uniqueness in law and the pathwise uniqueness for stochastic
differential equations,
{\it Theory Probab. Appl.} {\bf 48} (2003), no.3, 406--419.

\bibitem{dpz} 
Da Prato, G., Zabczyk, J.: 
{\it Stochastic Equations in Infinite Dimensions}, 
Encyclopedia of Mathematics and its Applications 44,
Cambridge University Press, 1992.

\bibitem{dpz2} 
Da Prato, G., Zabczyk, J.: 
{\it Ergodicity for infinite dimensional systems}, 
LMS Lecture Notes 229, Cambridge University Press, 1996.

\bibitem{cosa}
Ferrario, B.:
Absolute continuity of laws for semilinear stochastic equations with
additive noise, 
{\it Commun. on Stoch. Anal.} {\bf 2} (2008), no.2, 209--227;
Erratum, to appear in  
{\it Commun. on Stoch. Anal.} (2010).

\bibitem{gir}
Girsanov, I. V.:
On transforming a certain class of stochastic processes by absolutely
continuous substitution of measures, 
{\it Theory Probab. Appl.}  {\bf 5} (1960), no.3, 285--301.

\bibitem{IW}
Ikeda, N., Watanabe, S.:
{\it Stochastic differential equations and diffusion processes},
North-Holland Publishing Co., 1981.

\bibitem{KS}
Karatzas, I., Shreve, S. E.:
{\it Brownian motion and stochastic calculus},
Springer, 1988.

\bibitem{Ka}
Kazamaki, N.:
On a problem of Girsanov.  
{\it T\^ohoku Math. J.}  {\bf 29}  (1977), no.4, 597--600. 

\bibitem{K}
Krylov, N. V.:
A simple proof of a result of A. Novikov, 
e-print {\it arXiv:math/0207013v2} (2009). 


\bibitem{LS}
Liptser, R. S., Shiryaev, A. N.: 
{\it Statistics of random processes. I. General theory},
 Springer, 1977.

\bibitem{LS2}
Liptser, R. S., Shiryaev, A. N.: 
{\it Statistics of random processes. II. Applications},
 Springer, 1978.

\bibitem{nov}
Novikov, A. A.:
On an identity for stochastic integrals,
{\it Theory Probability Appl.}  {\bf 17} (1972), no.4, 717--720.


\end{thebibliography}

\end{document}